\documentclass{article}

\usepackage{amsfonts,amsmath,amssymb}
\usepackage[francais,english]{babel}
\usepackage[dvips]{graphicx}
\usepackage[T1]{fontenc}
\begin{document}
\newcommand{\und}{\vec{1}_d}
\newcommand{\suite}{_{n\geq 1}}
\newcommand{\ls}{\limsup}
\newcommand{\li}{\liminf}
\newcommand{\e}{\epsilon}
\newcommand{\eni}{\po \e_{n,i}\pf_{n\geq1, i \le m_n}}
\newcommand{\tab}{_{n\geq1,\;i\le m_n}}
\newcommand{\lb}{\newline}
\newcommand{\dd}{\delta}
\newcommand{\ld}{\llcorner}
\newcommand{\lr}{\lrcorner}
\newcommand{\alp}{\alpha}
\newcommand{\lab}{\lambda}
\newcommand{\Lab}{\Lambda}
\newcommand{\gam}{\gamma}
\newcommand{\Gam}{\Gamma}
\newcommand{\sig}{\sigma}
\newcommand{\Sig}{\Sigma}
\newcommand{\mmi}{\mid\mid}
\newcommand{\mmmi}{\mid\mid\mid}
\newcommand{\mmii}{\mid\mid_\infty}
\newcommand{\gMid}{\Bigg |}
\newcommand{\Mid}{\Big |}
\newcommand{\Mmi}{\Big|\Big|}
\newcommand{\Mmii}{\Big|\Big|_\infty}
\newcommand{\T}{\Theta}
\newcommand{\tht}{\theta}
\newcommand{\sli}{\sum\limits}
\newcommand{\sliin}{\sum\limits_{i=1}^n}
\newcommand{\sliid}{\sum\limits_{i=1}^d}
\newcommand{\sliik}{\sum\limits_{i=1}^k}
\newcommand{\sliiN}{\sum\limits_{i=1}^N}
\newcommand{\slijk}{\sum\limits_{j=1}^k}
\newcommand{\proliin}{\prod\limits_{i=1}^n}
\newcommand{\proliid}{\prod\limits_{i=1}^d}
\newcommand{\proliik}{\prod\limits_{i=1}^k}
\newcommand{\proliiN}{\prod\limits_{i=1}^N}
\newcommand{\prolijk}{\prod\limits_{j=1}^k}
\newcommand{\ili}{\int\limits}
\newcommand{\proli}{\prod\limits}
\newcommand{\limn}{\lim_{n\rightarrow\infty}\;}
\newcommand{\limk}{\lim_{k\rightarrow\infty}\;}
\newcommand{\inv}{\frac{1}}
\newcommand{\lsn}{\limsup_{n\rightarrow\infty}}
\newcommand{\lin}{\liminf_{n\rightarrow\infty}}
\newcommand{\lsk}{\limsup_{k\rightarrow\infty}}
\newcommand{\lik}{\liminf_{k\rightarrow\infty}}
\newcommand{\bculi}{\bigcup\limits}
\newcommand{\bcali}{\bigcap\limits}
\newcommand{\wt}{\widetilde}
\newcommand{\indep}{{\bot}\kern-0.9em{\bot}}
\newcommand{\beq}{\begin{equation}}
\newcommand{\Var}{\mathrm{Var}}
\newcommand{\eeq}{\end{equation}}
\newcommand{\mcal}{\mathcal}
\newcommand{\sq}{\sqrt}
\newcommand{\rar}{\rightarrow}
\newcommand{\lar}{\leftarrow}
\newcommand{\cvps}{\rightarrow_{p.s.}\;}
\newcommand{\cvpr}{\rightarrow_{P}\;}
\newcommand{\cvloi}{\rightarrow_{\mcal{L}}\;}
\newcommand{\ii}{\mathbf{i}}
\newcommand{\jj}{\vec{j}}
\newcommand{\kk}{\vec{k}}
\newcommand{\zzi}{z_{\mathbf{i}}}
\newcommand{\zzij}{z_{\mathbf{i},\vec{j}}}
\newcommand{\zzI}{\vec{z}_I}
\newcommand{\zzIJ}{\vec{z}_{I,J}}
\newcommand{\zz}{z}
\newcommand{\nk}{n_k}
\newcommand{\nkm}{n_{k-1}}
\newcommand{\nkk}{n_{k+1}}
\newcommand{\nkd}{n_k}
\newcommand{\nkkd}{n_{k+1}}
\newcommand{\Ln}{L_{n,1}}
\newcommand{\Lnn}{L_{n,2}}
\newcommand{\Lnl}{L_{n,\ell}}
\newcommand{\Pn}{\Pi_{n,1}}
\newcommand{\Pnn}{\Pi_{n,2}}
\newcommand{\Pnl}{\Pi_{n,\ell}}
\newcommand{\mn}{\vec{m}_n}
\newcommand{\Mn}{\vec{M}_n}
\newcommand{\hhn}{h_{n,1}}
\newcommand{\hhnn}{h_{n,2}}
\newcommand{\hhk}{h_{n_k,1}}
\newcommand{\hhkk}{h_{n_k,2}}
\newcommand{\mk}{\vec{m}_k}
\newcommand{\Mk}{\vec{M}_k}
\newcommand{\hhnl}{h_{n,\ell}}
\newcommand{\hnk}{h_{n_k}}
\newcommand{\bn}{b_{n,1}}
\newcommand{\bnn}{b_{n,2}}
\newcommand{\bnl}{b_{n,\ell}}
\newcommand{\bnk}{b_{n_k}}
\newcommand{\bnkk}{b_{n_{k+1}}}
\newcommand{\hkl}{{h_{n_k}^{(l)}}}
\newcommand{\hnl}{{h_n^{(l)}}}
\newcommand{\Tn}{T_{n,1}}
\newcommand{\Tnn}{T_{n,2}}
\newcommand{\Tnl}{T_{n,\ll}}
\newcommand{\hamn}{\hat{m}_n}
\newcommand{\harn}{\hat{r}_n}
\newcommand{\hafn}{\hat{f}_n}
\newcommand{\srl}{\stackrel}
\newcommand{\wap}{(\Omega,\mathcal{A},\rm I\kern-2pt P)}
\newcommand{\aoo}{\left\{}
\newcommand{\aff}{\right\}}
\newcommand{\coo}{\left[}
\newcommand{\cff}{\right]}
\newcommand{\poo}{\left(}
\newcommand{\pff}{\right)}
\newcommand{\po}{\left(}
\newcommand{\pf}{\right)}
\newcommand{\ao}{\left\{}
\newcommand{\af}{\right\}}
\newcommand{\co}{\left[}
\newcommand{\cf}{\right]}
\newcommand{\pooo}{\left(}
\newcommand{\pfff}{\right)}
\newcommand{\aooo}{\left\{}
\newcommand{\afff}{\right\}}
\newcommand{\cooo}{\left [}
\newcommand{\cfff}{\right]}
\newcommand{\poooo}{\left(}
\newcommand{\pffff}{\right)}
\newcommand{\aoooo}{\left\{}
\newcommand{\affff}{\right\}}
\newcommand{\coooo}{\left[}
\newcommand{\cffff}{\right]}
\newcommand{\FF}{\mathcal{F}}
\newcommand{\TT}{\mathcal{T}}
\newcommand{\GG}{\mathcal{G}}
\newcommand{\BBGG}{\mathcal{B}(\mathcal{G})}
\newcommand{\CC}{\mathcal{C}}
\newcommand{\KK}{\mathcal{K}}
\newcommand{\SSS}{\mathcal{S}}
\newcommand{\BB}{\mathcal{B}}
\newcommand{\PP}{\mathcal{P}}
\newcommand{\HH}{\mathcal{H}}
\newcommand{\NN}{\mathcal{N}}
\newcommand{\MM}{\mathcal{M}}
\newcommand{\DD}{\mathcal{D}}
\newcommand{\cc}{\widetilde{c}}
\newcommand{\EEE}{\mathbb{E}}
\newcommand{\NNN}{\mathbb{N}}
\newcommand{\PPP}{ \mathbb{P}}
\newcommand{\CCC}{\mathbb{C}}
\newcommand{\KKK}{\mathbb{K}}
\newcommand{\RRR}{\mathbb{R}}
\newcommand{\wtf}{\widetilde{f}}
\newcommand{\wth}{\mathfrak{h}}
\newcommand{\wtn}{\widetilde{n}}
\newcommand{\wtv}{\widetilde{v}}
\newcommand{\wtA}{\widetilde{A}}
\newcommand{\wtC}{\widetilde{C}}
\newcommand{\wtE}{\widetilde{E}}
\newcommand{\ovg}{\overline{g}}
\newcommand{\ovh}{\overline{h}}
\newcommand{\ovE}{\overline{E}}
\newcommand{\wtG}{\widetilde{G}}
\newcommand{\ovH}{\overline{H}}
\newcommand{\ovI}{\overline{I}}
\newcommand{\ovJ}{\overline{J}}
\newcommand{\ovK}{\overline{K}}
\newcommand{\wtN}{\widetilde{N}}
\newcommand{\wtP}{\widetilde{P}}
\newcommand{\ovR}{\overline{R}}
\newcommand{\wtPPP}{\widetilde{\mathbb{P}}}
\newcommand{\ovPPP}{\overline{\mathbb{P}}}
\newcommand{\wtF}{\widetilde{F}}
\newcommand{\wtI}{\widetilde{I}}
\newcommand{\wtK}{\widetilde{K}}
\newcommand{\wtFF}{\widetilde{\mathcal{F}}}
\newcommand{\ovFF}{\overline{\mathcal{F}}}
\newcommand{\cov}{\mathrm{Cov}}
\newcommand{\kif}{k\rightarrow\infty}
\newcommand{\nif}{n\rightarrow\infty}
\newcommand{\FFGG}{\FF\times\GG}
\newcommand{\vk}{\vskip10pt}
\newcommand{\Tproj}{\mathcal{T}_{0}}
\newcommand{\TTd}{\mathcal{T}_d}
\newcommand{\farc}{\frac}
\newcommand{\Nf}{\nabla_f}
\newcommand{\Nfn}{\nabla_f(\log_2n )}
\newcommand{\Cf}{\chi_f}
\newcommand{\nono}{\nonumber}
\newcommand{\Dn}{\Delta_n}
\newcommand{\DPn}{\Delta\Pi_n}
\newcommand{\norm}{\mid\mid \cdot \mid\mid}
\newcommand{\Xt}{(X(t))_{t\in T}}
\newcommand{\zklj}{z_{k,l,j}}
\newcommand{\znlj}{z_{n,l,j}}
\newcommand{\Idd}{[0,1)^d}
\newcommand{\pipi}{p_{i,n,\ii_0}}
\newcommand{\ppd}{2^{-pd}}
\newcommand{\iideux}{1\prec\mathbf{i}\prec 2^p}
\newcommand{\hk}{h_{n_k}}
\newcommand{\zik}{z_{i,n_{k}}}
\newtheorem{theo}{Theorem}
\newtheorem{ptheo}{Proof of theorem}
\newtheorem{lem}{Lemma}[section]
\newtheorem{plem}{Proof of lemma}[section]
\newtheorem{prop}{Property}[section]
\newtheorem{preuveprop}{Proof of property}[section]
\newtheorem{defi}{Definition}[section]
\newtheorem{propo}{Proposition}[section]
\newtheorem{popo}{Proof of proposition}[section]
\newtheorem{coro}{Corollary}[theo]
\newtheorem{pcoro}{Proof of corollary}[theo]
\newtheorem{rem}{Remark}[subsection]
\newtheorem{Fact}{Fact}[section] \numberwithin{equation}{section}

\title{A nonstandard uniform functional limit law for the increments of the multivariate empirical distribution function}
\author{Davit VARRON, Université de Franche-Comté}
\maketitle
\begin{abstract}
Let $(Z_i)_{i\geq 1}$ be an independent, identically distributed
sequence of random variables on $\RRR^d$. Under mild conditions on
the density of $Z_1$, we provide a nonstandard uniform functional
limit law for the following processes on $[0,1)^d$:
$$\Delta_n(z,h_n,\cdot):=s\mapsto \frac{\sliin
1_{[0,s_1]\times\ldots\times[0,s_d]}\poo\frac{Z_i-z}{h_n^{1/d}}\pff}{c\log
n},\;s\in [0,1)^d,$$ along a sequence $(h_n)\suite$ fulfilling
$h_n\downarrow 0,\;nh_n\uparrow,\;nh_n/\log c\rar c>0$. Here $z$
ranges through a compact set of $\RRR^d$. This result is an
extension of a theorem of Deheuvels and Mason \cite{DeheuvelsM2}
to the multivariate, non uniform case.\lb
\end{abstract}
\section{Introduction and statement of the result}
In this paper, we consider an independent, identically distributed
sequence of random vectors $(Z_i)_{i\geq 1}$ having a density $f$ on
an open set $O\subset \RRR^d$. We make the following assumption on
$f$:
\begin{tabbing}$ \hskip 5pt $ \= $ (H f) $ $ \hskip 5pt $
\= $f$ is continuous and strictly positive on $O$.
\end{tabbing}
Throughout this article, $s,s'\in \RRR^d$, we shall write $s\prec
s'$ when $s_i\le s'_i$ for each $i=1,\ldots,n$. Intervals and semi
intervals are implicitly understood as product of intervals or
semi intervals, namely \begin{align} \nono[s,s']:=&\{u\in
\RRR^d,\;s\prec u\prec
s'\}\\
=&[s_1,s'_1]\times\ldots\times[s_d,s'_d],\;s=(s_1,\ldots,s_d),\;s'=(s'_1,\ldots,s'_d).\end{align}
We shall also write $a\prec s$ (resp. $s\prec a$) for $s\in
\RRR^d$ and $a\in \RRR$ when $a\le s_i$ (resp. $s_i\le a$) for
each $i=1,\ldots,d$. For fixed $0<h<1$ and $z\in O$, we define the
following process on $[0,1)^d$:
$$
\mathbf{\Delta}_n(z,h,s):=\frac{1}{n}\sliin
1_{[0,s]}\poo\frac{Z_i-z}{h^{1/d}}\pff,\;s\in [0,1)^d.$$ These
processes, usually called functional increments of the empirical
distribution function, have been intensively investigated in the
literature (see, e.g., Shorack and Wellner \cite{MSW}, Van der
Vaart and Wellner \cite{Vander}, Deheuvels and Mason
\cite{DeheuvelsM2,DeheuvelsM3}, Einmahl and Mason
\cite{EinmahlM2}, Mason \cite{Mason2}). A particular domain of
investigation of these increments is when their almost sure
behavior is studied along a sequence of bandwidths $(h_n)\suite$
satisfying the following conditions:
\begin{tabbing}
$\hskip 20pt $ \= $ (HVE1)$ $\hskip 5pt $ \= $0<h_n<1,\;h_n\downarrow 0,\; nh_n\uparrow \infty, $\\
\> $ (HVE2)$ \> $nh_n/\log n\rar c$. \\
\end{tabbing}
Here, $c>0$ denotes a finite constant. Such conditions on the
sequence $(h_n)\suite$ are called Erdös-Rényi conditions, since
these two authors have given a pioneering result in this domain
(see \cite{ErdösRenyi}). Deheuvels and Mason \cite{DeheuvelsM2}
showed that, whenever the $(Z_i)_{i\geq 1}$ are uniformly
distributed on $[0,1]$, and under $(HVE1)-(HVE2)$, the increments
$n\mathbf{\Delta}_n(z,h,.)/(c\log n)$ have a nonstandard almost
sure behaviour. Before citing their result, we need to introduce
the following notations. Set $ B(\Idd)$ as the cone of all bounded
increasing functions $g$ on $\Idd$ (implicitly with respect to the
order $\prec$), satisfying $g(0)=0$. We shall endow this cone with
the topology spawned by the usual sup-norm $\mmi
g\mmi:=\sup_{s\in\Idd}\mid g(s)\mid$. Define the usually called
Chernoff function $h$ as
\beq h(x):=\left\{%
\begin{array}{ll}
    x\log x - x+1, & \hbox{for x>0;} \\
    1, & \hbox{for x=0;} \\
    \infty, & \hbox{for x<0.} \\
\end{array}%
\right.   \label{hChernoff}\eeq That function is known to play an
important role in the large deviation of Poisson processes on
$[0,1]$ (see, e.g., \cite{Lynch}). Define the following (rate)
function on $ B(\Idd)$. Whenever $g\in B(\Idd)$ is absolutely
continuous with respect to the Lebesgue measure on $\Idd$, we set
\beq I(g):= \ili_{\Idd} h(g'(s))ds\label{defI},\eeq $g'$ denoting
(a version of) the derivative of $g$ with respect to the Lebesgue
measure. Whenever $g$ fails to be absolutely continuous, we set
$I(g)=\infty$. Also define, for any $a>0$ ,  \beq \Gam_a:=\ao g\in
B(\Idd),\;I(g)\le 1/a\af.\label{Gama}\eeq In a pioneering work,
Deheuvels and Mason \cite{DeheuvelsM2} established the following
non standard uniform functional limit law for the
$\mathrm{\Delta}_n(z,h_n,\cdot)$, when the $(Z_i)$ are uniform on
$[0,1]$.
\begin{theo}[Deheuvels, Mason, 1992]
Assume that $d=1$ and that the $(Z_i)_{i\geq 1}$ are uniformly
distributed on $[0,1]$. Let $0\le a<b<1$ be two real numbers, and
let $(h_n)\suite$ be a sequence of positive constants satisfying
$(HVE1)-(HVE2)$ for some constant $c>0$. Then we have almost surely
\begin{align}
\nono&\limn \sup_{z\in [0,1-h_n]}\inf_{g\in \Gam_c} \Mmi
\frac{n}{c\log n}\mathbf{\Delta}_n(z,h_n,\cdot)-g\Mmi=0,\\
\nono\forall g\in \Gam_c,\;&\limn \inf_{z\in [0,1-h_n]}\Mmi
\frac{n}{c\log n}\mathbf{\Delta}_n(z,h_n,\cdot)-g\Mmi=0.
\end{align}
\end{theo}
As a corollary, the authors showed that, when the sequence of
bandwidth $(h_n)\suite$ satisfies $(HVE1)-(HVE2)$, the
Parzen-Rosenblatt kernel density estimator is \textbf{not}
uniformly strongly consistent. They proved this non-consistency
result by making use of some optimisation techniques on Orlicz
balls (see Deheuvels and Mason \cite{DeheuvelsM4}). The aim of the
present paper is to provide a generalisation of the former result
to the case where the $(Z_i)_{i\geq 1}$ take values in $\RRR^d$.
This generalisation can be stated as follows.
\begin{theo}\label{T3}
Assume that the $(Z_i)_{i\geq 1}$ have a density $f$ satisfying $(H
f)$. Let $H\subset O$ be a compact set with nonempty interior. Let
$(h_n)\suite$ be a sequence of positive constants fulfilling
$(HVE1)$ and $(HVE2)$. Then we have almost surely
\begin{align}
&(i)\;\forall z\in H, \forall g \in \Gam_{cf(z)},\limn \inf\aoo
\Mmi\Dn(z',h_n,\cdot)-g\Mmi ,\; z'\in H\aff=0,\\
&(ii)\limn\sup_{z\in H}\;\inf\aoo \Mmi\Dn(z,h_n,\cdot)-g\Mmi,\;
g\in\Gam_{cf(z)}\aff=0.
\end{align}
\end{theo}Denote by $f_n(K,z,h_n)$ the usual kernel density
estimator with bandwidth $h_n$ and kernel $K$. A consequence of
Theorem \ref{T3} is that, under $(HVE1)-(HVE2)$, $f_n(K,z,h_n)$ is
not uniformly consistent (in a strong sense) over (say) an
hypercube of $\RRR^d$.\lb \textbf{Corollary}: Let $K$ be a kernel
with compact support and bounded variation. Assume $(Hf)$ and
$(HVE1)-(HVE2)$. Let $H\subset O$ be a compact with nonempty
interior. Then the following event holds with probability one:
$$\exists \e>0,\;\exists n_0,\;\forall n\geq n_0,\;\sup_{z\in H}
\mid f_n(K,z,h_n)-f(z)\mid >\e.$$ \textbf{Proof}: The proof
follows exactly the lines of Deheuvels and Mason (see
\cite{DeheuvelsM2}, Theorem 4.2) and is based on some optimisation
results on Orlicz Balls that have been provided in Deheuvels and
Mason \cite{DeheuvelsM4}. $\Box$\lb From now on, we shall make use
of the following notation
$$\Delta_n(z,h_n,s):=\;\frac{\sliin
1_{[0,s]}\poo\frac{Z_i-z}{h_n^{1/d}}\pff}{cf(z)\log n},\;s\in
\Idd.\label{Dn}$$
\begin{rem}\end{rem}
Deheuvels and Mason \cite{DeheuvelsM6} have already given a
nonstandard functional limit law for a single increment
$\Dn(z_0,h_n,\cdot)$ when $(HVE 2)$ is replaced by $nh_n/\log\log
n\rar c>0$. Their result is presented in a more general setting,
considering the $\Delta_n(z_0,h_n,\cdot)$ as random measures
indexed by a class of sets.\lb The remainder of this paper is
organised as follows. In \S \ref{deux} we provide some tools in
large deviation theory, which are consequences of results of
Arcones \cite{Arcones1} and Lynch and Sethuraman \cite{Lynch}. In
\S \ref{trois}, a uniform large deviation principle for
"poissonized" versions of the $\Delta_n(z,h_n,\cdot)$ is
established. In \S \ref{quatre} and \S \ref{cinq}, we make use of
the just-mentioned uniform large deviation principle to prove
Theorem \ref{T3}.
\section{Uniform large deviation principles}\label{deux} The main tool
we shall make use of in \S \ref{quatre} and  \S\ref{cinq} is a
uniform large deviation principle for a triangular array of
compound Poisson processes. We must first remind some usual
notions in large deviation theory. Let $(E,d)$ be a metric space.
A real function $J: E\rar [0,\infty]$ is said to be a rate
function (implicitly for $(E,d)$) when the sets $\{ x\in
E:\;J(x)\le a\},\;a\geq 0$, are compact sets of $(E,d)$. We shall
first show that $I$ is a rate function on $\po B(\Idd),\norm\pf$
by approximating it by suitably chosen simple rate functions.
\subsection{Approximations of $I$}
Given $g\in B(\Idd)$ and a Borel set $A$, we shall write \beq
g(A):=\ili_{\Idd}1_{A}d g,\label{mesure}\eeq which is valid as
soon as either $g$ or $1_A$ has bounded variation. For any integer
$p\geq 1$ and for each $1\prec \ii\prec 2^p$ set \beq A^p_{\ii}:=
2^{-p}\left[\ii-1,\ii\right),\label{Aiip}\eeq with the notation
$\ii-1:=(i_1-1,\ldots,i_d-1)$. Recall that $h$ is given in
(\ref{hChernoff}), and that $\lab$ is the Lebesgue measure on
$\Idd$. The following functions will play the role of
approximations of $I$ (given in (\ref{defI})), as $p\rar \infty$ :
\begin{align} I_p(g):=&\sli_{1\prec \ii \prec 2^p}2^{-pd}h\po
2^{pd}g( A^p_{\ii})\pf\label{Ip}\\
\nono=&\sli_{1\prec \ii \prec 2^p}\lab\po A^p_{\ii}\pf
h\poo\frac{g( A^p_{\ii})}{\lab( A^p_{\ii})}\pff,\;g\in B(\Idd).
\end{align}
We point out the following properties of the function $I$.
\begin{propo}\label{PropI}
For each $g\in  B(\Idd)$, we have \beq\lim_{p\rar \infty}
I_p(g)=I(g).\label{etcacom}\eeq Moreover, $I$ is a rate function
on $\po B(\Idd),\norm\pf$.
\end{propo}
\addtocounter{popo}{1}\textbf{Proof}: Choose $g\in B(\Idd)$
arbitrarily and assume that $I(g)>0$ (nontrivial case). In a first
time, we suppose that $g$ has bounded variation, so that it can be
interpreted as a finite measure. Denote by $\TT_p$ the
$\sig$-algebra of $\Idd$ spawned by the sets $A^p_{\ii},\;1\prec
\ii\prec 2^p$. Clearly, for all $p\geq 1$, the measure $g$ is
absolutely continuous with respect to the (trace of the) Lebesgue
measure $\lab$ on $\TT_p$. Furthermore, the corresponding
Radon-Nicodym derivative is given by the following equality. \beq
L_p:=\frac{d g}{d\lab}\mid_{\TT_p} =\sli_{1\prec \ii \prec 2^p}
1_{A^p_{\ii}}\frac{g(A^p_{\ii})} {\lab(A^p_{\ii})}\label{deroc}.
\eeq Clearly the $\sig$-algebra spawned by the (increasing)
sequence $(\TT_p)_{p\geq 1}$ is equal to the Borel $\sig$-algebra
of $\Idd$. Assume first that $g$ is absolutely continuous with
respect to $\lab$. According to Dacunha-Castelle and Duflo
\cite{Dacu}, p. 63, the sequence $L_p$ converges $\lab+g$ almost
everywhere to a positive function $L$ satisfying $L=g'\;$($\lab+g$
almost everywhere). Now select $0<l<I(g)$ arbitrarily. By
definition of $I$, there exists $\e>0$ satisfying
\[\ili_{\e<L<1/\e}h(L)d\lab>l.\]
Since $L_p\rar L$ ($\lab+g$ almost everywhere as $p\rar \infty$) and
since $h$ is continuous, we have
\[\liminf_{p\rar\infty}h(L_p)1_{\{\e<L_p<1/\e\}}\geq
h(L)1_{\{\e<L<1/\e\}}\;\ \lab+g\;almost\;everywhere\] Hence by an
application of Fatou's lemma,
\[\liminf_{p\rar\infty}
\ili_{\e<L_p<1/\e}h(L_p)d\lab\geq\ili_{\e<L<1/\e}h(L)d\lab>l.\]
Since $\sup_{p\geq 1} I_p(g)\le I(g)$ by a straightforward use of
Jensen's inequality, and since $l<I(g)$ was chosen arbitrarily, we
readily infer that $I_p(g)\rar I(g)$ as $p\rar \infty$. Now assume
that $I(g)=\infty$ and that $g$ is not absolutely continuous with
respect to $\lab$. According to Dacunha-Castelle and Duflo
\cite{Dacu}, p. 63, the sequence $L_p$ converges $\lab+g$ almost
everywhere to a positive function $L$ satisfying $(\lab+g)( \{
L=\infty\})=:\tau>0$. Define
\[\ell(x):=x^{-1}h(x)=\log (x)-1+x^{-1} ,\;x>0.\]
Clearly, $ \ell(x)\rar\infty$ as $\mid x\mid \rar \infty$. Now
select $l>0$ arbitrarily, and choose $A>0$ satisfying
\[\inf_{x>A}\ell(x)>\frac{2l}{\tau}.\]
Since $L_p\rar L$ ($\lab+g$ almost everywhere as $p\rar \infty$) we
have $g(L_p>A)>\tau/2$ for all large $p$, whence
\begin{align} \nono I_p(g)  \geq& \ili_{L_p\in
(A,\infty)}\ell(L_p)L_p\;d\lab\\
\nono = & \ili_{L_p\in (A,\infty)}\ell(L_p)d g\\
\nono\geq& \frac{2l}{\tau}g(L_p>A)\\
 >&l\label{guito}.
\end{align}
We have shown that (\ref{etcacom}) is true for each $g$ with
bounded variation. Whenever $g$ has infinite variation, then it
can be shown that $I_p(g)\rar \infty$ by a discrete version of the
argument that have just been invoked to obtain (\ref{guito}). We
omit details for sake of briefness.\lb Since all the functions
$I_p$ are $\norm$-continuous and since $I_p(g)\uparrow I(g)$ for
all $g\in B(\Idd)$, we conclude that $I$ is lower-semicontinuous
for $\norm$. Hence, $I$ is a rate function if and only if the set
$\Gam_a$ is totally bounded for each $a>0$ (recall (\ref{Gama})).
Since $x^{-1}h(x)\rar \infty$ as $\mid x\mid \rar \infty$, we
have, for some constant $M>0$, \beq \mid x \mid \le \mid x \mid
1_{\mid x\mid \le M}+h(x),\eeq from where we readily infer that
\beq \ili_{\Idd} \mid g'\mid d\lab\le M+1/a\text{ for each
}a>0\text{ and }g\in \Gam_a.\eeq
 Applying the Arzela-Ascoli criterion,
we conclude that, for each $a>0$, the closed set $\Gam_a$ is
totally bounded, which entails that $I$ is a rate function on $\po
B(\Idd),\norm\pf$. This concludes the proof of Proposition
\ref{PropI}.$\Box$ \subsection{Uniform large deviations in $\po
B(\Idd),\norm\pf$} We shall now give a definition of a large
uniform large deviation principle in the metric space $\po
B(\Idd), \norm\pf$. In the sequel, $(\e_{n,i})_{n\geq 1,i\le m_n}$
will always denote a triangular array of positive numbers
satisfying $\max_{i\le m_n} \e_{n,i}\rar 0$ as $\nif$. Let
$(X_{n,i})\tab$ be a triangular array of random elements on
probability space $\po\omega,\TT',\PPP\pf$, taking values in
$B(\Idd)$. In order to handle carefully the notions of inner and
outer probabilities, we shall that each $X_{n,i}$ is a suitable
projection mapping from $\po\Omega,\TT'\pf$ to $E$, where
$$\Omega:=
\proli_{n=1}^{\infty}\proli_{i=1}^p
 B(\Idd),\;\;\TT':=\bigotimes\limits_{n=1}^\infty\bigotimes\limits_{i=1}^p
\TT,$$  and $\TT$ is the Borel $\sig$-algebra of $\po
B(\Idd),\norm\pf$. From now on, outer and inner probabilities
$\PPP^*$ and $\PPP_*$ are understood with $(\Omega,\TT')$ as the
underlying probability space. We say that $(X_{n,i})\tab$
satisfies the Uniform Large Deviation Principle (ULDP) for
$(\e_{n,i})\tab$ and for a rate function $J$ whenever the two
following conditions hold.
\begin{itemize}
\item For any $\norm$-open set $O\subset  B(\Idd)$ we have
\beq\lin \min_{i\le m_n} \e_{n,i}\log\poo\PPP_*\po
X_{n,i}(\cdot)\in O\pf\pff\geq -J(O)\label{GDU2}.\eeq \item For
any $\norm$-closed set $F\subset B(\Idd)$ we have \beq\lsn
\max_{i\le m_n}\e_{n,i}\log\poo\PPP^*\po X_{n,i}(\cdot)\in
F\pf\pff\le -J(F)\label{GDU1}.\eeq
\end{itemize}
\begin{rem}\end{rem} The same definition holds for triangular arrays of random variables taking
values in $\RRR^p,\;p\geq 1$. The norm $\norm$ can then be
replaced by any norm.\vk Arcones \cite{Arcones1} provided a
powerful tool to establish Large Deviation Principles for
sequences of bounded stochastic processes. Some verifications lead
to the conclusion that the just-mentioned tool can be used in our
context. Recall that the sets $A_{\ii}^p$ have been define by
(\ref{Aiip}). Consider the following finite grid, for $p\geq 1$ :
\beq s_{\ii,p}:=2^{-p}(\ii-1),\;1\prec \ii\prec 2^p.\eeq Given,
$p\geq 1$ and $g\in B(\Idd)$, we write
$$g^{(p)}= \sli_{1\prec \ii \prec 2^p}1_{A_{\ii}^p}g(s_{\ii,p}).$$
\begin{propo}\label{propo: ecrit}
Let $(X_{n,i})\tab$ be a triangular array of random elements
taking values in $(B(\Idd))$ almost surely, and let
$(\e_{n,i})\tab$ be a triangular array of positive real numbers.
Assume that the following conditions are satisfied.
\begin{enumerate}
\item   The
 triangular array of stochastic process $( X_{n,i}^{(p)})\tab$
satisfies the ULDP for $(\e_{n,i})\tab$ and for the rate function
$I_p$ on $\po  B(\Idd),\norm\pf$. \item For each $\tau>0$ and
$M>0$ there exists $p\geq 1$ satisfying
\[\lsn \max_{i\le m_n}\e_{n,i}\log\poo\PPP^*\poo\max_{1\prec \ii \prec 2^p}\sup_{s\in A_{\ii}^p}\mid X_{n,i}(t)-X_{n,i}(s_{\ii}^p)\mid \geq
\tau\pff\pff\le -M.\]
\end{enumerate}
Then $(X_{n,i})\tab$ satisfies the ULDP for $(\e_{n,i})\tab$ and for
the following rate function.
$$J(g):=\sup_{p\geq
1}I_p\poo g^p\pff,\;g\in  B(\Idd).$$
\end{propo}
\textbf{Proof}: The proof follows exactly the same lines as in the
proof of Theorem 3.1 of Arcones \cite{Arcones1}. Using theses
arguments in our context remains possible since the cone $
B(\Idd)$ is a closed subset of $L^{\infty}(\Idd)$ for the usual
sup norm $\norm$. We avoid writing the proof for sake of
briefness. $\Box$\lb Another tool we shall make an intensive use
of is a ULDP for random vectors with mutually independent
coordinates.
\begin{propo}\label{propo: espaceproduit}
Let $(X_{n,i})_{n\geq 1,\;1\le i\le m_n}$  and $(Y_{n,i})_{n\geq
1,\;1\le i\le m_n}$ be two triangular arrays of random vectors
taking values in $\RRR^d$ and $\RRR^{d'}$ respectively, and
satisfying $X_{n,i}\;\indep \;Y_{n,i}$ for each $n\geq 1,\;1\le
i\le m_n$. Assume that both $(X_{n,i})_{n\geq 1,\;1\le i\le m_n}$
and $(Y_{n,i})_{n\geq 1,\;1\le i\le m_n}$ satisfy the ULDP for a
triangular array $(\e_{n,i})\tab$ and for two rate functions $J_1$
and $J_2$ respectively. Then the triangular array
$(X_{n,i},Y_{n,i})_{n\geq 1, i\le m_n}$ satisfies the ULDP for
$(\e_{i,n})\tab$ and for the following rate function.
\[J(z_1,z_2):=J_1(z_1)+J_2(z_2),\; z_1\in \RRR^d, \;z_2\in
\RRR^{d'}.\]
\end{propo}
\textbf{Proof}: The proof follows the same lines as Lemma 2.6 and
Corollary 2.9 in Lynch and Sethuraman \cite{Lynch}. In the
just-mentioned article, the authors make use of the notions of
Weak Large Deviation Principle and of LD-tightness for sequences
of random variables in a Polish space. These notions can be easily
extended to the frame of triangular arrays of random variables.
$\Box$\lb The following proposition is nothing else than the
contraction principle in the framework of ULDP (see, e.g.,
\cite{Arcones1}, Theorem 2.1 for the most general version of that
principle).
\begin{propo}\label{contraction}
Let $(X_{n,i})\tab$ be a triangular arrays of $\RRR^p$ valued
random vectors satisfying the ULDP for a triangular array
$(\e_{n,i})\tab$ and for a rate function $J$. Let $\mcal{R}$ be a
continuous mapping from $\RRR^d$ to $\po  B(\Idd),\norm\pf$. Then
$(\mcal{R}(X_{n,i}))\tab$ satisfies the ULDP for $(\e_{n,i})\tab$
and for the following rate function.
$$J_{\mcal{R}}(g):=\inf\{J(x),\;\mcal{R}(x)=g\},\;g\in
 B(\Idd),$$ with the convention $\inf \emptyset =\infty$.\lb
\end{propo}
\textbf{Proof}: Straightforward. $\Box$\lb
 The
following proposition shall be useful in our the proof of our
Lemma \ref{ingsl}.
\begin{propo}\label{Chernoffunivorme}
Let $ (X_{n,i})_{n\geq 1,i\le m_n}$ be a triangular array of real
random variables and let $(\e_{n,i})_{n\geq 1,i\le m_n}$ be a
triangular array of positive real numbers. Assume that there
exists a strictly convex positive function $J$ on $\RRR$ and a
real number $\mathbf{\mu}$ such that $J(\mathbf{\mu})=0$ and
\begin{align}
\forall a>\mathbf{\mu},&\;\limn \max_{i\le m_n}
\Mid\e_{n,i}\log\poo\PPP\poo
X_{n,i}\geq a\pff\pff-J(a)\Mid=0,\label{dede1}\\
\forall a<\mathbf{\mu},&\;\limn \max_{i\le m_n}
\Mid\e_{n,i}\log\poo\PPP\poo X_{n,i}\le
a\pff\pff-J(a)\Mid=0\label{dede2}.
\end{align}
Then $\po X_{n,i}\pf\tab$ satisfies the ULDP for
$\po\e_{n,i}\pf\tab$ and for $J$.
\end{propo}
\textbf{Proof}: The proof is routine calculus.$\Box$
\section{A ULDP for poissonised versions of the $\Delta_n(z,h_n,\cdot)$}\label{trois}
Define the following process, for each integer $n\geq 1$. \beq
\DPn(z,h_n,s):=\frac{\sli_{i=1}^{\eta_n}
1_{[0,s]}\poo\frac{Z_i-z}{h_n^{1/d}}\pff}{c f(z)\log
n}\label{DPn},\;s\in \Idd.\eeq Here $\eta_n$ is a Poisson random
variable independent of $(Z_i)_{i\geq 1}$, with expectation $n$.
These "poissonized" versions of the processes $\Dn(z,h_n,\cdot)$
can be identified to random (Poisson) measures by the following
relation \beq \DPn(z,h_n,A):=\ili_{\Idd} 1_A(s)
d\DPn(z,h_n,s),\;A\text{ Borel}\label{DPnA}.\eeq The key of our
proof of Theorem \ref{T3} is the following ULDP.
\begin{propo}\label{gdpommes} Let $(z_{i,n})_{n\geq 1,\;1\le i\le
m_n}$ be a triangular array of elements of $H$. Under the
assumptions of Theorem \ref{T3}, the triangular array of processes
$(\DPn(z_{i,n},h_n,\cdot))_{n\geq 1,\;1\le i\le m_n}$ satisfies
the ULDP in $( B(\Idd),\mmi \cdot \mmi)$ for the rate function $I$
and for the following triangular array \beq
\e_{n,i}:=\frac{1}{cf(z_{i,n})\log n},\;n\geq 1,\;1\le i\le
m_n.\label{epsic}\eeq
\end{propo}
\begin{rem}\end{rem} Proposition \ref{gdpommes} is true whatever the constant $c>0$ appearing in assumption
(HVE1). This remark will show up to be useful in Lemma
\ref{encore2} in \S \ref{cinq}.\lb
\addtocounter{popo}{1}\textbf{Proof}: To prove proposition
\ref{gdpommes}, we shall make use of Proposition \ref{propo:
ecrit}. We hence have to check conditions 1, 2 and 3 of the
just-mentioned proposition. This will be achieved through several
lemmas.
\subsection{A preliminary lemma}
Recall notation (\ref{mesure}). To check condition 2 of
Proposition \ref{propo: ecrit}, we need first to establish the
following lemma.
\begin{lem}\label{ingsl}Assume that the hypothesis of Theorem
\ref{T3} are satisfied. Then, for each $p\geq 1$ and for each
$1\prec \ii_0\prec 2^p$, the triangular array of random variables
$(\DPn(z_{i,n},h_n,A_{\ii_0}^p))_{n\geq 1,\;1\le i\le m_n}$
 satisfies the ULDP in $[0,\infty)$ for the triangular array $(\e_{n,i})_{n\geq 1,\:i\le m_n}$
 and for the following rate function:
\beq \wtI_p(x):=2^{-pd}h\po\frac{x}{2^{-pd}}\pf=\lab\po
A_{\ii_0}^p\pf h\poo \frac{x}{\lab\po
A_{\ii_0}^p\pf}\pff,\;x\geq0.\label{wtIp}\eeq
\end{lem}
\addtocounter{plem}{1}\textbf{Proof}: Fix once for all $p\geq 1$
and $1\prec \ii_0\prec 2^d$. We shall make use of Proposition
\ref{Chernoffunivorme}, with $J:=\wtI_p$ and $\mu:=2^{-pd}$. We
give details only for the proof of (\ref{dede1}), as proving
(\ref{dede2}) is very similar. Fix $a>2^{-pd}$. For each integers
$n\geq 1$ and $1\le i\le m_n$, we set (recall (\ref{DPnA}))
\begin{align}
\nono V_{i,n,\ii_0}:=&cf(z_{i,n})(\log
n)\DPn(z_{i,n},h_n,A_{\ii_0}^p),\\
\nono p_{i,n,\ii_0}:=&\PPP\po Z_1\in
z_{i,n}+h_n^{1/d}A_{\ii_0}^p\pf\label{pinii0}.
\end{align}
Clearly $V_{i,n,\ii_0}$ is a Poisson random variable with
expectation $n\pipi$. Since the density $f$ satisfies $(H f)$ and
since $\lab(A_{\ii_0}^p)=\ppd$, we have \beq \limn \max_{1\le i\le
m_n} \Mid \frac{p_{i,n,\ii_0}}{f(z_{i,n})\ppd h_n}-1
\Mid=0.\label{zuzu}\eeq Hence according to (HVE2) we have,
ultimately as $\nif$, \beq\min_{1\le i\le m_n}\frac{a c
f(z_{i,n})\log n}{np_{i,n,\ii_0}}>1.\label{igt1}\eeq We then make
use of Chernoff's inequality for Poisson random variables to get,
for all large $n$ (satisfying (\ref{igt1})) and  for all $1\le
i\le m_n$,
\begin{align}
\nono \PPP\poo \DPn(z_{i,n},h_n,A_{\ii_0}^p)\geq a\pff=\;&\PPP\poo
V_{i,n,\ii_0}\geq acf(z_{i,n})\log n\pff\\
 \le\;& \exp\poo -np_{i,n,\ii_0}h\poo\frac{acf(z_{i,n})\log
n}{np_{i,n,\ii_0}}\pff\pff.\label{igt2}
\end{align}
But (\ref{igt2}) in combination with (\ref{zuzu}) entails \beq
\lsn \max_{1\le i\le m_n}\frac{p_{i,n,\ii_0}}{f(z_{i,n})
h_n}h\poo\frac{acf(z_{i,n})\log n}{np_{i,n,\ii_0}}\pff\le \ppd
h\poo\frac{a}{\ppd}\pff\label{igt3},\eeq which, together with
(\ref{igt2}) leads to \beq\lsn \max_{1\le i\le
m_n}\e_{n,i}\log\poo\PPP\poo\DPn(z_{i,n},h_n,A_{\ii_0}^p)\geq
a\pff\pff\le -\wtI_p(a).\label{babadou1}\eeq Now select $y>a$
arbitrarily. If we could show that
\[\lin \min_{1\le i\le m_n} \e_{n,i}\log\poo\PPP\poo \DPn(z_{i,n},h_n,A_{\ii_0}^p)\geq a\pff\pff\geq
-\wtI_p(y),\]
 then, as $y>a$ was chosen arbitrarily, and since  $\wtI_p$ is increasing on $[a,\infty)$,
 we should be able to conclude the proof of (\ref{dede1})with $J=\wtI_p$. Now
set $\phi(t):=\exp\po\exp(t)-1\pf,\;t\in \RRR$ and notice that
$h(z)=\max_{u\in \RRR} zu -\log \po\phi (u)\pf$ for each $z>0$. Set
$u_0:=\log(2^{pd} y)$, so as
\begin{align}
h(2^{pd}y)
=&2^{pd}yu_0-\log\po\phi(u_0)\pf.\label{zuzu2}\end{align} Denote
by $F$ the distribution function of a Poisson random variable with
expectation 1, and define $F_0$ by  \beq dF_0(x):={\phi\po
u_0\pf}^{-1}\exp(u_0x) d F(x).\label{F0}\eeq Let "*" be the
convolution operator for infinitely divisible laws and notice
that, for each $L>0$, we have \begin{align}
dF_0^{*L}(\cdot)=&{\phi\po
u_0\pf}^{-L}\exp(u_0\cdot)dF^{*L}(\cdot)\label{ouioui1},\\
\EEE_{F_0^{*L}}(X)=&2^{pd}Ly\label{espoir},\\
\Var_{F_0^{*L}}(X)=&L\Var_{F_0}(X)\label{variance}
\end{align}
Here we have written $\EEE_{F}(X)$ as the expectation of a random
variable with distribution $F$. Now fix $\dd>0$ satisfying
$[y-\dd,y+\dd]\subset [a,\infty[$ arbitrarily. Obviously,
$F^{*n\pipi}$ is the distribution function of $cf(z_{i,n})(\log
n)\;\DPn(z_{i,n},h_n,A_{\ii_0}^p)$, whence
\begin{align}
\nono &\PPP\poo \DPn(z_{i,n},h_n,A_{\ii_0}^p)\geq a\pff\\
\nono\geq\;&\PPP\poo \DPn(z_{i,n},h_n,A_{\ii_0}^p)\in[y-\dd,y+\dd]\pff\\
\nono= \;&\ili_{\frac{x}{c f(z_{i,n})\log n}\in [y-\dd,y+\dd]}dF^{*n\pipi}(x) \\
\nono\geq \;& \exp\poo-u_0(y+\dd) cf(z_{i,n})\log n\pff\times\ili_{\frac{x}{c f(z_{i,n})\log n}\in [y-\dd,y+\dd]}\exp(u_0x)dF^{*n\pipi}(x)\\
\nono\geq\;& \exp\poo-  cf(z_{i,n})(\log n)u_0(y+\dd)+np_{i,n,\ii_0}\log\po\phi( u_0)\pf\pff\\
&\;\times\ili_{\frac{x}{c f(z_{i,n})\log n}\in [y-\dd,y+\dd]}dF_0^{*np_{i,n,\ii_0}}(x)\label{pacs2}\\
\nono:=&\;a_{i,n,\ii_0,\dd}\times b_{i,n,\ii_0,\dd}.
\end{align}
Here (\ref{pacs2}) is a consequence of (\ref{ouioui1}), with
$L:=n\pipi$. Now let $n\geq 1$ be an integer large enough to
fulfill (recall (\ref{zuzu}))  \beq\max_{1\le i\le
m_n}\Mid\frac{n\pipi}{\ppd cf(z_{i,n})\log n}-1\Mid \le
u_0\log\po\phi(u_0)\pf^{-1}\dd\label{qop1},\eeq which enables us
to write the following chain of inequalities.
\begin{align}
\nono &c f(z_{i,n})(\log n) u_0(y+\dd)-n\pipi \log\po\phi(u_0)\pf
\\
\nono\le&\;
\ppd (y+\dd)cf(z_{i,n})\log n\;\poo u_02^{pd}-\log\po\phi(u_0)\pf+ u_0\dd \pff\\
\nono\le& \;\ppd cf(z_{i,n})\log n \poo h\po2^{pd}y\pf+
u_0(2^{pd}+1)\dd\pff\\
\nono=&\; cf(z_{i,n}) \log n\poo \wtI_p(y)+\ppd \po2^{pd}+1\pf u_0\dd\pff\\
\le\; &cf(z_{i,n})\log n \poo \wtI_p(y)+2 u_0\dd\pff.
\end{align}
Therefore we have, for all large $n$ and for all $1\le i \le m_n$,
\beq a_{i,n,\ii_0,\dd}\geq \exp\poo -cf(z_{i,n})\log n
\po\wtI_p(y)+2u_0\dd\pf\pff\label{despa1},\eeq where
$u_0=\log(2^{pd}y)$ depends on $y>a$ only. It remains to show that
 \beq \limn
 \min_{1\le i\le m_n} b_{i,n,\ii_0,\dd}=1.\label{despa2}\eeq
Consider $n$ large enough to fulfill (recall (\ref{zuzu}))
\beq\nono \frac{y-\dd}{y+\ppd \dd}<\min_{1\le i\le m_n}
\frac{n\pipi}{\ppd c f(z_{i,n})\log n}\le\max_{1\le i\le m_n}
\frac{n\pipi}{\ppd c f(z_{i,n})\log n}< \frac{y+\dd}{y-\ppd
\dd},\eeq so as, for all $1\le i \le m_n$, \beq \frac{n\pipi}{\ppd
c f(z_{i,n})\log n}\times [y-\ppd\dd,y+\ppd \dd]\subset\;
]y-\dd,y+\dd[\label{grouve1},\eeq and hence
\begin{align}
\nono b_{i,n\ii_0,\dd}\geq &\;
\ili_{\frac{x}{n\pipi}\in[2^{pd}y-\dd,2^{pd}y+\dd]}dF_0^{*n\pipi}(x).
\end{align}
Recalling (\ref{espoir}) and (\ref{variance}) we get, by the
Bienaymé-Tchebychev inequality,
\begin{align}
\label{fert1} 1-b_{i,n,\ii_0,\dd}\le \frac{\Var_{F_0}(X)}{\dd
n\pipi}.
 \end{align}
By assumption $(H f)$ we infer that the $h_n^{-1}\pipi$ are
bounded away from zero, from where (\ref{pacs2}) follows. Then
(\ref{pacs2}), (\ref{despa1}) and (\ref{despa2}) entail \beq\lin
\min_{1\le i\le m_n}\e_{n,i}\log \poo \PPP\poo
\DPn(z_{i,n},h_n,A_{\ii_0}^p)\geq a\pff\pff\geq
-\wtI_p(y)-2u_0\dd\label{balor}.\eeq  Assertion (\ref{dede1}) is
then proved by combining (\ref{babadou1}) with (\ref{balor}), as
$\dd>0$ is arbitrary. $\Box$\vk
\subsection{Verification of condition 2 of Proposition \ref{propo: ecrit}}
For $n\geq 1$ and $1\le i\le m_n$, define the following
$\RRR^{2^{pd}}$ valued random vector:
\begin{align}\nono X_{n,i}:=&\poo X_{\ii_0,n,i}\pff_{1\prec \ii_0\prec 2^p}\\
\nono :=& \poo \DPn(z_{i,n},h_n,A_{\ii_0}^p)\pff_{1\prec
\ii_0\prec 2^p}. \end{align} Notice that the random variables
$X_{i_0,n,i},\;1\prec \ii_0\prec 2^p$ are mutually independent for
fixed $n\geq 1$ and $1\le i \le m_n$ by usual properties of
Poisson random measures. Hence, by Lemma \ref{ingsl} together with
Proposition \ref{propo: espaceproduit} we deduce that the
triangular array $(X_{n,i})\tab$ satisfies the ULDP with
$(\e_{n,i})\tab$ and with the following rate function. \beq
I'_p(x):=\sli_{\iideux}\ppd h\poo\frac{x_{\ii}}{\ppd}\pff,\;x\in
{[0,\infty)}^{2^{pd}}.\label{aeift}\eeq Here we have written
$x:=(x_{\ii})_{1\prec \ii\prec 2^p}$. We now define the following
mappings from ${[0,\infty)}^{2^{pd}}$ to $\poo
 B(\Idd)\pff$
\[\begin{array}{rcl} \mcal{R}_p(x):\Idd&\mapsto& [0,\infty)\\
s&\rar&\sli_{A_{\ii}^p\subset[0,s]}x_{\ii}.
\end{array}\]
Denote by $[x]$ the integer part of a real number $x$ ($[x]\le
x<[x]+1$), and write $[s]:=([s_1],\ldots,[s_d])$
 for any $s=(s_1,\ldots,s_d)\in \RRR^d$. We point out that with probability one (recall the notations of Proposition \ref{propo: ecrit})
\begin{align}
\nono \mcal{R}_p(X_{n,i})(s)=& \DPn\poo z_{i,n},h_n,
2^{-p}[2^ps]\pff\\
\nono= &{\DPn\poo z_{i,n},h_n,s\pff}^{(p)},\;s\in \Idd.
\end{align}
For fixed $p \geq 1$, we make use of the contraction principle
(Proposition \ref{contraction}) to conclude that
$(\mcal{R}_p(X_{n,i}))\tab$ satisfies the ULDP for
$(\e_{n,i})\tab$ and for the following rate function.  \beq
\ovI_p(g):=\inf\aoo
I'_p(x),\;x\in{[0,\infty)}^{2^{pd}},\;\mcal{R}_p(x)=g\aff,\;g\in
 B(\Idd)\label{resto1},\eeq with the convention $\inf \emptyset
=\infty$. Obviously, the set appearing in (\ref{resto1}) is non
void if and only if $g$ is the cumulative distribution function of
a purely atomic measure with atoms belonging to the grid
$\{s_{\ii,p},\;1\prec \ii\prec 2^p\}.$ In that case we have
\[ \ovI_p(g)=\sli_{1\prec \ii \prec 2^p}\ppd h\poo\frac{g(A_{\ii}^p)}{\ppd}\pff=I_p(g).\]
Here, we have identified $g$ to a positive finite measure on
$\Idd$ (recall (\ref{mesure})). Assumption 2 of Proposition
\ref{propo: ecrit} is then satisfied.
\subsection{Verification of condition 3 of Proposition \ref{propo:
ecrit}} Fix $\tau>0$ and $M>0$. We have to prove that, provided
that $p$ is large enough,
\begin{align}
\nono &\lsn\max_{1\le i\le
m_n}\e_{n,i} \\
\nono&\log\poo \PPP\poo\max_{1\prec \ii\prec 2^p}\sup_{s\in
A_{\ii}^p} \Mid\DPn\po z_{i,n},h_n,s\pf
-\DPn\po z_{i,n},h_n,2^{-p}(\ii-1)\pf \Mid\geq \tau\pff\pff\\
\le& -M\label{judok}.\end{align} For fixed $p\geq 1,\;n\geq 1,\;1\le
i\le m_n$, a rough upper bound gives
\begin{align}
\nono&\PPP\poo\max_{1\prec \ii\prec 2^p}\sup_{s\in A_{\ii}^p} \Mid\DPn\po
z_{i,n},h_n,s\pf -\DPn\po z_{i,n},h_n,2^{-p}(\ii-1)\pf
\Mid\geq \tau\pff\\
  \nono\le& 2^{pd}\max_{1\prec \ii\prec 2^p}\PPP\poo\mathop{\sup_{2^{-p}(\ii-1)\prec s}}_{\prec 2^{-p}\ii} \Mid\DPn\po
z_{i,n},h_n,s\pf -\DPn\po z_{i,n},h_n,2^{-p}(\ii-1)\pf \Mid\geq
\tau\pff\\
\nono\le&\;\PPP\poo\DPn\po z_{i,n},h_n,2^{-p}\ii\pf-\DPn\po
z_{i,n},h_n,2^{-p}(\ii-1)\pf\geq
\tau\pff\\
 =:&\;\PPP_{i,n,\ii,p}\label{judok4}.
\end{align}
We shall now write
\begin{align}
\nono W_{i,n,\ii,p}:=&cf(z_{i,n})\log n\pooo\DPn\po
z_{i,n},h_n,2^{-p}\ii\pf-\DPn\po
z_{i,n},h_n,2^{-p}(\ii-1)\pf\pfff,\\
\nono\mu_{i,n,\ii,p}:=&\PPP\poo\frac{Z_1-z_{i,n}}{h_n^{1/d}}\in
[0,2^{-p}\ii)-[0,2^{-p}(\ii-))\pff,\text{ and}\\
\nu_{\ii,p}:=&\lab\poo [0,2^{-p}\ii)-[0,2^{-p}(\ii-))\pff\le
d2^{-p}.\label{abala}\end{align} Clearly, $W_{i,n,\ii,p}$ is a
Poisson random variable with expectation $n\mu_{i,n\ii,p}$.
Moreover, by assumption $(H f)$ we have \beq \limn
\mathop{\min_{1\le i\le m_n,}}_{\iideux}\frac{cf(z_{i,n})(\log
n)\nu_{\ii,p}}{n\mu_{i,n,\ii,p}}=1\label{jb1}.\eeq Recall that
$x^{-1}h(x)\rar\infty$ as $x\rar\infty$. We can then choose
$A_{M,\tau}>1$ large enough to satisfy \beq \inf_{x\geq
A_{M,\tau}}\frac{h(x)}{x}>\frac{8M}{\tau}.\label{rafz1}\eeq By
(\ref{abala}) we can choose $p$ large enough to fulfill \beq
\min_{\iideux}
\frac{\tau}{2\nu_{\ii,p}}>A_{\tau,M}.\label{rafz3}\eeq Assertion
(\ref{jb1}) together with (\ref{rafz3}) leads to the following
inequality, for all large $n$, for all $1\le i\le m_n$ and for all
$\iideux$.
 \beq \frac{cf(z_{i,n})\tau\log
n}{n\mu_{i,n,\ii,p}}\geq
\frac{\tau}{2\nu_{\ii,p}}>A_{\tau,M}>1.\label{lpmz}\eeq Applying
Chernoff's inequality to the Poisson random variables
$W_{i,n,\ii,p}$ we get, for all large $n$ and for all $1\le i\le
m_n$,
\begin{align}
\nono\PPP_{i,n,\ii,p}=&\;\PPP\poo W_{i,n,\ii,p}\geq
\tau c f(z_{i,n})\log n\pff\\
\nono\le& \exp \pooo-n\mu_{i,n,\ii,p}h\poo \frac{cf(z_{i,n})\tau\log
n}{n\mu_{i,n,\ii,p}}\pff\pfff.
\end{align}
Therefore, recalling (\ref{jb1}) and (\ref{lpmz}), the following
inequality holds for all large $n$, for all $1\le i\le m_n$ and
for all $1\prec \ii\prec 2^p$. \begin{align} \nono
\PPP_{i,n,\ii,p}\le&
\exp\poo-\frac{1}{2}cf(z_{i,n})\nu_{\ii,p}(\log
n) h\po\frac{\tau}{2\nu_{\ii,p}}\pf\pff\\
\le&\exp\po-cf(z_{i,n})2M\log n\pf.\label{rafz4}
\end{align}
Here, (\ref{rafz4}) is a consequence of (\ref{rafz3}). By combining
(\ref{rafz4}) with and (\ref{judok4}) we get, for all
large $n$ and for each $1\le i\le m_n$,
\begin{align}
\nono&\PPP\poo\max_{1\prec \ii\prec 2^p}\sup_{s\in A_{\ii}^p}
\Mid\DPn\po z_{i,n},h_n,s\pf -\DPn\po z_{i,n},h_n,2^{-p}(\ii-1)\pf
\Mid\geq \tau\pff\\
\nono\le &\exp\po-2Mcf(z_{i,n})\log n+\log(2^{pd})\pf,\end{align}
which proves (\ref{judok}) and shows that condition 3 of
Proposition \ref{propo: ecrit} is satisfied, as $f$ is bounded
away from zero on $H$. We can now make use of the just-mentioned
proposition in combination with Proposition \ref{PropI} to
conclude the proof of Proposition \ref{gdpommes}. $\Box$
\section{Proof of part (i) of Theorem \ref{T3}}\label{quatre}
Denote by $\mathrm{Int}(H)$ the interior of $H$, and fix $z\in
\mathrm{Int}(H)$, $g\in \Gam_{cf(z)}$, and $\e>0$. We set \beq
g^{\e}:=\aoo g'\in  B(\Idd),\;\mmi g'-g \mmi <
\e\aff.\label{geps}\eeq By lower semi continuity  of $I$ in $\po
 B(\Idd),\norm\pf$ (recall Proposition \ref{PropI}), there exists
$\alp_1>0$ satisfying \beq I\po
g^\e\pf=\frac{1-3\alp_1}{cf(z)}.\label{re6}\eeq Now choose an
hypercube with nonempty interior
$H':=[a_1,b_1]\times\ldots\times[a_p,b_p]$ fulfilling $H'\subset
H$, $\PPP\poo Z_1\in H'\pff\le 1/2$ and \beq\inf_{z'\in
H'}\frac{f(z')}{f(z)}
> \frac{1-2\alp_1}{1-\alp_1}\label{re7}.\eeq Such a choice is possible since $H$ has a nonempty interior by assumption. We now divide
$H'$ into disjoint hypercubes $z_{i,n}+h_n^{1/d}{[0,1)}^d,\;1\le
i\le m_n$, where $m_n$ is the maximal number of disjoint
hypercubes we can construct without violating \beq
\bculi_{i=1}^{m_n}\aoo z_{i,n}+h_n^{1/d}{[0,1)}^d\aff\subset H'
\label{dicret31}.\eeq Notice that, as $n\rar \infty$, \beq
m_n=h_n^{-1+o(1)}=n^{(1+o(1))}.\label{re8}\eeq Now recall
$(\ref{DPn})$. By making use of a well-known "poissonization"
technique (see, e.g., Mason \cite{Mason1}, Fact 6), we get the
following upper bound for all large $n$.
\begin{align} \nono& \PPP\poo \bcali_{z'\in H}\aoo
\Dn(z',h_n,\cdot)\notin
g^{\e}\aff\pff\\
\nono
\le\;&\PPP\poo\bcali_{i=1}^{m_n}\aoo\Dn(z_{i,n},h_n,\cdot)\notin
g^{\e}\aff\pff \\
 \le\;& 2\PPP\poo\bcali_{i=1}^{m_n}\aoo\DPn(z_{i,n},h_n,\cdot)\notin
g^{\e}\aff\pff\label{palot2}\\
 =&\;2\proli_{i=1}^{m_n}\poo1-\PPP\poo\DPn(z_{i,n}h_n,\cdot)\in
g^{\e}\pff\pff\label{palot1}\\
\le\;& 2\exp\poo-m_n\min_{1\le i\le
m_n}\PPP\poo\DPn(z_{i,n},h_n,\cdot)\in g^{\e}\pff\pff\label{gast}
\end{align}
The transition between (\ref{palot2}) and (\ref{palot1}) is a
classical property of Poisson random measures, while inequality
(\ref{gast}) is a consequence of $1-u\le \exp(-u),\;u\geq 0$. We
now make use of Proposition \ref{gdpommes} (with the open ball
$g^{\e}$) to get, for all large $n$ (recall (\ref{re6})),
\begin{align}
\nono\PPP\poo \bcali_{z'\in H}\aoo \Dn(z',h_n,\cdot)\notin
g^{\e}\aff\pff\le&2\exp \poo-m_n\;\min_{1\le i\le
m_n}\;n^{-\frac{f(z_{i,n})}{f(z)}(1-2\alp_1)}\pff\\
\nono \le & \exp \po -n^{\alp_1}\pf,\end{align} which is a
consequence of (\ref{re7}) and (\ref{re8}). Hence we conclude by
the Borel-Cantelli lemma that, almost surely,
\[\limn \inf\aoo
\mmi\Dn(z',h_n,\cdot)-g\mmi ,\; z'\in H\aff\le \e.\] As $\e>0$ was
chosen arbitrarily, the proof of part (i) of Theorem \ref{T3} is
concluded for each $z\in \mathrm{Int}(H)$. Now the case where
$z\in H$ does not belong to $\mathrm{Int}(H)$ is treated by making
use of the following argument: for each $z_1\in H,\;g_1\in
\Gam_{cf(z_1)}$ and $\e>0$, there exists $z_2\in \mathrm{Int}(H)$
and $g_2\in \Gam_{cf(z_2)}$ satisfying $\mmi g_1-g_2\mmi <\e$.
Such an argument is valid by $(H f)$ and by Lemma \ref{titiegro}
(see below).$\Box$
\section{Proof of part (ii) of Theorem \ref{T3}}\label{cinq}
We shall make use of somewhat usual blocking arguments along the
following subsequence $n_k:=\co\exp(k/\log k)\cf,\;k \geq 3$ and
its associated blocks $N_k:=\{n_{k-1}+1,\ldots,n_k\}$. Given
$A\subset B(\Idd)$ and $\e>0$ we shall write \beq A^\e:=\aoo g\in
 B(\Idd),\;\inf_{g'\in A}\mmi
 g-g'\mmi<\e\aff\label{Aeps2}.\eeq
The following lemma shall come in handy.
\begin{lem}\label{titiegro}
For any $\e>0$ and $L>0$ there exists $\eta>0$ satisfying, for each,
$L'\in [(1+\eta)^{-1}L,L]$, $\Gam_{L'}\subset \Gam_L^{\e}.$
\end{lem}
\textbf{Proof}: The proof is routine analysis.$\Box$\lb Now fix
$\e>0$. Since  $I$ is lower-semi continuous on
$\po B(\Idd),\norm\pf$ (recall Proposition \ref{PropI}) we deduce
that, given $z\in H$, there exists $\alp_z>0$ satisfying \beq
I\poo B(\Idd)- \Gam_{cf(z)}^\e\pff=\frac{1+3\alp_z}{cf(z)}
\label{hoi1}.\eeq By $(H f)$ and Lemma  \ref{titiegro} we can
construct an hypercube $H_z$ with nonempty interior satisfying the
following
conditions. \begin{align}&z\in H_z,\;\;H_z\subset O,\label{boup} \\
&\inf_{z_1,z_2\in H_z}\frac{f(z_1)}{f(z_2)}\geq
\frac{1+\alp_z}{1+2\alp_z}\label{hoi2},\\& \bculi_{z'\in
H_{z}}\Gam_{cf(z')}\subset \Gam_{cf(z)}^\e,\label{hoi8}\\&\PPP\poo
Z_1\in \bculi_{z\in H_{z}}\aoo z+ {[0,{\hk}^{1/d})}^d \aff\pff\le
1/2.\label{hoii}\end{align} The compact set $H$ is included in the
union of the interiors of $H_z,\;z\in H$, from where we can
extract a finite union, noted as \beq H\subset\bculi_{l=1}^{L}
\mathrm{Int}{H}_{z_l}\subset\bculi_{l=1}^{L} H_{z_l}\subset
O.\label{gfdsm}\eeq Our problem is now reduced to showing that,
for fixed $l=1,\ldots,L$,
 \beq \lsn \sup_{z\in H_{z_l}}\inf_{g\in
\Gam_{cf(z_l)}} \mmi \Dn(z,h_n,\cdot)-g\mmi\le
10\e\label{hoi3}\;\;almost\;surely.\eeq We now fix $1\le l\le L$,
and we write $H_{z_l}=:[a_1,b_1]\times\ldots\times[a_d,b_d]$. We
now introduce a parameter $\dd>0$ that will be chosen in function
of $\e$ in the sequel. For each $k\geq 1$, we cover $H_{z_l}$ by
hypercubes \beq H_{z_l}\subset\bculi_{1\le i\le
m_{n_k}}C_{i,n_k}\subset O,\label{discret3},\eeq with
\begin{align} \nono C_{i,n_k}:=&\zik+{[0,(\dd
h_{n_k})^{1/d})}^d,\;k\geq 1,\;1\le i\le
m_{n_k}\text{ and}\\
 m_{n_k}:=&\proli_{p=1}^d \poo\coo \frac{b_p-a_p}{(\dd
\hk)^{1/d}}\cff+1\pff.\label{re1}\end{align} Now define, for each
$k\geq 1$, $n\in N_k$, $z\in H$,
$$\mcal{H}_n(z,s):=\frac{1}{c\log n_k}\sli_{i=1}^n 1_{[0,s)}\poo\frac{Z_i-z}{{h_{n_k}}^{1/d}}\pff,\;s\in \Idd.$$
We shall first show that, for any choice $\dd>0$, we have almost
surely \beq\lsn \sup_{1\le i\le m_{n_k}}\inf_{g\in
\Gam_{cf(z_l)}}\mmi \mcal{H}_n(z_{i,n_k},\cdot)-g\mmi\le
2\e\label{re3}.\eeq Consider the following probabilities for all
large $k$.
\[\PPP_k:=\PPP\poo\bculi_{1\le i\le m_{n_k}}\bculi_{n\in N_k}
\mcal{H}_n(z_{i,n_k},\cdot)\notin \Gam_{cf(z_l)}^{2\e}\pff.\] We
have, ultimately as $\kif$,
\begin{align}\PPP_k\le m_k\max_{1\le i\le m_{n_k}}\PPP\poo\bculi_{n\in N_k}
\mcal{H}_n(z_{i,n_k},\cdot)\notin
\Gam_{cf(z_l)}^{\e}\label{hoi4}\pff.
\end{align}
We now make use of a well-known maximal inequality (see, e.g.,
Deheuvels and Mason \cite{DeheuvelsM2}, Lemma 3.4) to get, for all
large $k$ and for all $1\le i\le m_{\nk}$, \beq
\PPP\poo\bculi_{n\in N_k} \mcal{H}_n(z_{i,n_k},\cdot)\notin
\Gam_{cf(z_l)}^{2\e}\pff\le 2\PPP\poo
\mcal{H}_{n_k}(z_{i,n_k},\cdot)\notin
\Gam_{cf(z_l)}^{\e}\pff.\label{pars1}\eeq We point out that the
conditions of Lemma 3.4 in \cite{DeheuvelsM2} are satisfied since,
by a straightforward use of Markov's inequality we have,
ultimately as $\kif$,
\[\sup_{z\in H} \max_{n\in N_k}\PPP\poo
\mmi \mcal{H}_{n_k}(z,\cdot)-\mcal{H}_n(z,\cdot)\mmi\geq \e\pff\le
\frac{1}{2}.\] Making use of (\ref{pars1}) in (\ref{hoi4}), we
obtain, for all large $k$,
\begin{align} \nono \PPP_k\le&\; 2m_k\max_{1\le i\le
m_{n_k}}\PPP\poo \mcal{H}_{n_k}(z_{i,n_k},\cdot)\notin
\Gam_{cf(z_l)}^{\e}\pff\\
\nono=&\;2m_{n_k}\max_{1\le i\le m_{n_k}}\PPP\poo
\Delta_{n_k}(\zik,\hk,\cdot)\notin\Gam_{cf(z_l)}^{\e}\pff\\
\le &\;4m_{n_k}\max_{1\le i\le m_{n_k}}\PPP\poo
\Delta\Pi_{n_k}(\zik,\hk,\cdot)\notin\Gam_{cf(z_l)}^{\e}\pff.
\end{align}
The last inequality is a consequence of usual poissonization
techniques (see, e.g., Mason \cite{Mason1}, Fact 6). We now make
use of Proposition \ref{gdpommes}, which, together with
(\ref{hoi1}) leads to the following inequality, ultimately as
$\kif$,
 \[\PPP_k\le 4m_{n_k} \max_{1\le m_k}\exp\poo-\frac{f(z_{i,n_k})}{f(z_l)}(1+2\alp_{z_l})\log
 n_k\pff.\]
 Moreover (\ref{hoi2}) entails $\PPP_k\le 4m_{n_k}\exp\po-(1+\alp_{z_l})\log n_k\pf.$
Since $m_{n_k}=\hnk^{-1+o(1)}=n_k^{1+o(1)}$ as $k \rar \infty$
(recall (\ref{re1})), the sumability of $\PPP_k$ follows,
 which proves (\ref{re3}) by the Borel-Cantelli lemma. We point out that (\ref{re3}) is true whatever the choice of $\dd>0$ (recall (\ref{discret3})). We now focus on showing that, for a small value of
 $\dd>0$ we have
 \beq \lsk\; \sup_{z\in H_{z_l}}\min_{1\le i\le m_{n_k}} \max_{n\in
N_k}\mmi \mcal{H}_n(\zik,\cdot)-\Dn(z,h_n,\cdot)\mmi\le
7\e\;\;a.s,\label{relais1}\eeq which will be achieved through two
separate lemmas.
\begin{lem}\label{encore2}
Assume that the conditions of Theorem \ref{T3} are fulfilled.
There exists $\dd_\e>0$ such that, for any choice of
$0<\dd<\dd_\e$  we have almost surely \beq \nono\lsk \max_{n\in
N_k}\max_{1\le i\le m_{n_k}}\;\;\sup_{z\in C_{i,n_k}}\Mmi
\mcal{H}_n(\zik,\cdot)-\frac{f(z)}{f(\zik)}\mcal{H}_n(z,\cdot)\Mmi
\le\e.\eeq
\end{lem}
\textbf{Proof}: For all large $k$ we have
\begin{align}\nono &\PPP\poo\max_{n\in N_k}\max_{1\le i\le
m_{n_k}}\;\;\sup_{z\in C_{i,n_k}}\Mmi
\mcal{H}_n(\zik,\cdot)-\frac{f(z)}{f(\zik)}\mcal{H}_n(z,\cdot)\Mmi > \e\pff\\
\nono = &\PPP\poo \bculi_{1\le i\le m_{n_k}}\bculi_{n\in
N_k}\;\;\sup_{z\in C_{i,n_k}}\Mmi
\mcal{H}_n(\zik,\cdot)-\frac{f(z)}{f(\zik)}\mcal{H}_n(z,\cdot)\Mmi > \e\pff\\
 \le &m_{n_k}\max_{i\le m_{n_k}}\PPP\poo\bculi_{n\in
N_k}\;\;\sup_{z\in C_{i,n_k}}\Mmi
\mcal{H}_n(\zik,\cdot)-\frac{f(z)}{f(\zik)}\mcal{H}_n(z,\cdot)\Mmi >
\e\pff\label{rdde2}
\end{align}
Fix $k\geq 1$, $1\le i\le m_{n_k}$ and $z\in
\zik+(\dd\hk)^{1/d}\Idd$. We write $\zik:=(\zik^1,\ldots,\zik^d)$,
$z:=(z^1,\ldots,z^d)$ and $Z_j:=(Z_j^1,\ldots,Z_j^d),\;j\geq 1$.
Notice that for each $p=1,...,d$ we have $\zik^p\le z^p\le
\zik^p+(\dd\hk)^{1/d}$. Hence, in virtue of the equality $\mid
1_A-1_B\mid =1_{A-B}+1_{B-A}$ we have, for each integer $j$ we
have almost surely, for each $(s_1,\ldots,s_d)\in \Idd$,
\begin{align}
\nono&\Mid
1_{[0,s)}\poo\frac{Z_j-z}{\hk^{1/d}}\pff-1_{[0,s) }\poo\frac{Z_j-\zik}{\hk^{1/d}}\pff\Mid\\
\nono=&1_{\aoo \co z,z+\hk^{1/d}s\pf -\co \zik,\zik+\hk^{1/d}s\pf\aff}(Z_j)+1_{\aoo\co \zik,\zik+\hk^{1/d}s\pf- \co z,z+\hk^{1/d}s\pf\aff}(Z_j)\\
\nono\le &\sli_{l=1}^d1_{\co
\zik^l+s_l\hk^{1/d},\zik^l+\hk^{1/d}(s_l+\dd^{1/d})\cf
}(Z_j^l)\proli_{1\le
p\not = l\le d}1_{\co \zik^p,\zik^p+\hk^{1/d}(s_p+\dd^{1/d})\cf }(Z_j^p)\\
 &\;+\sli_{l=1}^d1_{\co \zik^l,\zik^l+(\dd\hk)^{1/d}\cf }(Z_j^l)\proli_{1\le p\not = l\le d}1_{\co \zik^p,\zik^p+\hk^{1/d}s_p\cf }(Z_j^p)\label{hero}\\
=&:X_{j,k,i,\dd}(s)\label{aub}.
\end{align}
Here (\ref{hero}) follows from $\zik^l\le z^l\le
\zik^l+\dd^{1/d}\hk^{1/d},\;l=1,\ldots,d$. As the
$X_{j,k,i,\dd}(\cdot)$ are positive processes almost surely,
(\ref{aub}) entails, for all large $k$ and for all $1\le i\le
m_{n_k}$,
\begin{align}
\nono &\PPP\Big{(}\bculi_{n\in N_k}\sup_{z\in C_{i,n_k}}\Mmi \mcal{H}_n(\zik,\cdot)-\frac{f(z)}{f(\zik)}\mcal{H}_n(z,\cdot)\Mmi > \e\Big{)}\\
\nono \le& \PPP\Big{(}\bculi_{n\in N_k}\sup_{z\in C_{i,n_k}}\; \sup_{s\in \Idd}\\
\nono &\;\;\;\;\;\;\sli_{j=1}^{n}\Mid 1_{[0,s)}\poo\frac{Z_j-z}{\hk^{1/d}}\pff-1_{[0,s)}\poo\frac{Z_j-\zik}{\hk^{1/d}}\pff\Mid\geq \e cf(\zik)\log n_k\Big{)}\\
\nono \le &\PPP\Big{(}\bculi_{n=1}^{\nk}\sup_{s\in\Idd} \sli_{j=1}^nX_{j,k,i,\dd}(s)\geq \epsilon
cf(\zik)\log n_k\Big{)}\\
\le & \PPP\poo\Mmi \sli_{j=1}^{\nk}X_{j,k,i,\dd}(\cdot)\Mmi\geq \e cf(\zik)\log n_k\pff\label{auba1}.
\end{align}
But a close look at (\ref{hero}) leads to the conclusion that,
almost surely, for each $s\in \Idd$,
\begin{align}
\nono 0\le& \;\sli_{j=1}^{n_k} X_{j,k,i,\dd}(s)\\
\le &\; 2d cf(\zik)\log n_k\mathop{\sup_{s,s'\in [0,2)^d,}}_{\mmi s'-s \mmi^d<\dd}\Mid\Delta_{n_k}(\zik,h_{n_k},s')-\Delta_{n_k}(\zik,h_{n_k},s)\Mid\label{ara}
\end{align}
Here we have written $\mid s\mid_d:=\max\{\mid
s_j\mid,\;j=1,\ldots,p\}$. Now (\ref{ara}) together with
(\ref{auba1}) entails
\begin{align}
\nono &\frac{1}{2}\PPP\poo\bculi_{n\in N_k}\sup_{z\in C_{i,n_k}}\Mmi
\mcal{H}_n(\zik,\cdot)-\frac{f(z)}{f(\zik)}\mcal{H}_n(z,\cdot)\Mmi > \e\pff\\
\nono \le& \frac{1}{2}\PPP\poo \mathop{\sup_{s,s'\in
[0,2)^d,}}_{\mmi s'-s
\mmi_d<\dd^{1/d}}\Mid\Delta_{n_k}(\zik,h_{n_k},s')-\Delta_{n_k}(\zik,h_{n_k},s)\Mid>\frac{\e}
{2 d}\pff\\
\le & \PPP\poo \mathop{\sup_{s,s'\in [0,1)^d,}}_{\mmi s'-s
\mmi<\frac{\dd^{1/d}}{2}}\Mid\Delta\Pi_{n_k}(\zik,h_{n_k},2s')-\Delta\Pi_{n_k}(\zik,h_{n_k},2s)\Mid>\frac{\e}
{2 d}\pff\label{concon2}
\end{align}
Here (\ref{concon2}) follows from poissonization techniques. Now
consider the following sequence $\wth_n:=2^{d}h_n,\;n\geq 1$.
Clearly, $(\wth_n)\suite$ satisfies (HVE1) and (HVE2), replacing
$c$ by $\mathfrak{c}:=2^dc$. Moreover, for each $k\geq 1$, $1\le
i\le m_{n_k}$ we have almost surely, for all $s\in \Idd$, \beq
\Delta\Pi_{n_k}(\zik,h_{n_k},2s)=\Delta\Pi_{n_k}(\zik,\mathfrak{h}_{\nk},s).\eeq
Applying Proposition \ref{gdpommes} we deduce that the triangular
array of processes
\[U_{k,i}(\cdot):=\Delta\Pi_{n_k}(\zik,\mathfrak{h}_{n_k},2\cdot),\;k\geq 1,\;1\le
i\le m_{n_k}\] satisfies the ULDP in $ ( B(\Idd),\norm)$ (see \S
\ref{deux}) for the rate function $I$ and for the following
triangular array:
\[\e_{k,i}:=(c2^df(\zik)\log n_k)^{-1}k\geq 1,\;1\le
i\le m_{n_k}.\]Now consider the following set
\[\Gam:=\aoo g\in \MM([0,1)^d),\;I(g)\le \frac{4}{2^dc\beta }\aff.\]
By proposition \ref{PropI}, there exists $\dd_{\e}>0$ such that
\beq \sup_{g\in 2^d\Gam}\;\sup_{s,s'\in[0,2)^d,\mmi s'-s\mmi_d\le
\dd_\e^d/2}\mid g(s')-g(s)\mid< (4d )^{-1}\e.\label{concon1}\eeq
Now choose $0<\dd<\dd_{\e}$ arbitrarily for the construction of
the $\zik,\;k\geq 1,1\le i\le m_{n_k}$ (recall (\ref{discret3})).
By lower-semicontinuity of $I$, the closed set
\[F:=\ao g\in\MM([0,2)^d),\; \inf_{g'\in \Gam}\mmi g-g'\mmi_{[0,2)^d}\geq \frac{2^{-d}\e}{8 d}\af\]
satisfies $I(F)>4/(2^dc\beta).$ Hence, (\ref{concon2}) together
with (\ref{concon1}) leads to the following inequalities for all
large $k$ and for each $1\le i\le m_{n_k}$.
\begin{align}
\nono &\PPP\poo\bculi_{n\in N_k}\sup_{z\in C_{i,n_k}}\Mmi
\mcal{H}_n(\zik,\cdot)-\frac{f(z)}{f(\zik)}\mcal{H}_n(z,\cdot)\Mmi > \e\pff\\
\nono \le & 2\PPP\poo \mathop{\sup_{s,s'\in [0,1)^d,}}_{\mmi s'-s
\mmi_d<\dd^{1/d}/2}\Mid U_{k,i}(s')-U_{k,i}(s)\Mid>\e
(2 d)^{-1}\pff\\
\nono \le& 2\PPP\poo\Delta\Pi_{n_k}(\zik,\mathfrak{h}_{n_k},\cdot)\in F\pff\\
\nono \le &2\exp \poo - \frac{3}{4}I\po F\pf
\mathfrak{c}f(\zik)\log
n_k\pff\\
\nono \le &2\exp\poo - 3\times \frac{c2^df(\zik)}{\beta
c2^d}\log n_k\pff\\
\le & 2\exp \po -3\log n_k\pf.\label{rdde}
\end{align} Now (\ref{rdde}) in combination with (\ref{rdde2}) entails, for all large $k$, \beq\PPP\poo\max_{n\in
N_k}\max_{i\le m_{n_k}}\sup_{z\in C_{i,n_k}}\Mmi
\mcal{H}_n(\zik,\cdot)-\frac{f(z)}{f(\zik)}\mcal{H}_n(z,\cdot)\Mmi
> \e\pff\le \frac{2m_{n_k}}{{n_k}^{3}}\label{souri}\eeq But for
fixed $\dd>0$ we have $m_{n_k}=\hk^{-1+o(1)}={n_k}^{1+o(1)}$ as
$\kif$. The proof of Lemma \ref{encore2} is concluded by applying
the Borel-Cantelli lemma to (\ref{souri}). $\Box$\vk
\begin{lem}\label{encore1}
Under the assumptions of Theorem \ref{T3}, for any choice of
$\dd>0$, we have almost surely
\[\lsk \max_{1\le i\le m_{n_k}}\sup_{z\in C_{i,n_k}}\;\max_{n\in
N_k}\Mmi\Delta_n(z,h_n,\cdot)-\frac{f(z)}{f(\zik)}\mcal{H}_n(z,\cdot)\Mmi\le
6\e.\]
\end{lem}
\textbf{Proof}: For all large $k$ and for all $1\le i\le m_{n_k}$,
$z\in C_{i,n_k}$, $n\in N_k$ we have almost surely, for each $s\in
\Idd$,
\beq\Delta_n(z,h_n,s)=T_{n,i,k}\frac{f(z)}{f(\zik)}\mcal{H}_n\poo
z,\rho_{n,k}s\pff,\label{rata}\eeq with $T_{i,n,k}:f(\zik)\log
n_k/f(z)\log n$ and $\rho_{n,k}^d:=\hk/h_n$. First notice that\beq
\limk\nono \max_{1\le i\le m_{n_k}}\;\sup_{z\in C_{i,n_k}}\mid
T_{n,i,k}-1\mid =0,\;\;\ \limk \max_{n\in N_k}\mid
\rho_{n,k}-1\mid=0.\label{topil}\eeq Moreover, by Proposition
\ref{PropI} we have \beq \lim_{T\rar 1,\rho\rar
1}\sup_{g\in\Gam_{cf(z_l)}}\mmi T
g(\rho^{1/d}\cdot)-g(\cdot)\mmi=0.\label{cgta}\eeq Finally, by
(\ref{re3}) and by Lemma \ref{encore2} we have, for all large $k$
and for all $1\le i\le m_{n_k}$, $z\in C_{i,n_k}$, $n\in N_k$,
\beq\inf_{g\in\Gam_{cf(z_l)}}\Mmi
\frac{f(z)}{f(\zik)}\mcal{H}_n(z,\cdot)-g\Mmi<
3\e\;\;almost\;surely.\label{funk}\eeq Hence, combining
(\ref{rata}), (\ref{topil}), (\ref{cgta}), (\ref{funk}) and the
triangle inequality, we obtain almost surely, for all large $k$ and
for all $n\in N_k$ :
\begin{align}
\nono&\Mmi\Delta_n(z,h_n,\cdot)-\frac{f(z)}{f(\zik)}\mcal{H}_n(z,\cdot)\Mmi\\
\nono \le& 6\e,
\end{align}
which proves Lemma \ref{encore1}. $\Box$ \vk \textbf{End of the
proof of part(ii) of Theorem \ref{T3}}: By combining Lemma
\ref{encore1} with Lemma \ref{encore2} we conclude that
(\ref{relais1}) is true for $\dd>0$ small enough. Now
(\ref{relais1}) together with (\ref{re3}) leads to
\[\lsn\sup_{z\in H_{z_l}}\inf_{g\in \Gam_{c f(z_l)}} \mmi
\Dn(z,h_n,\cdot)-g\mmi \le9\e\;\;almost\; surely.\] Whence,
recalling (\ref{hoi8}), \beq\lsn\sup_{z\in H_{z_l}}\inf_{g\in
\Gam_{c f(z)}} \mmi \Dn(z,h_n,\cdot)-g\mmi \le
10\e\;\;almost\;surely.\label{hzlpa}\eeq Repeating (\ref{hzlpa})
for each $l=1,\ldots, L$ (recall (\ref{gfdsm})) we get
\[\lsn\sup_{z\in H}\inf_{g\in \Gam_{c f(z)}} \mmi
\Dn(z,h_n,\cdot)-g\mmi \le 10\e\;\;almost\;surely.\] As $\e>0$ was
chosen arbitrarily, the proof of part(ii) of Theorem \ref{T3} is
concluded.$\Box$

\end{document}